\documentclass[preprint,aps]{revtex4-1}
\usepackage{amssymb,amscd,amsmath,amsthm,amsfonts}
\usepackage{latexsym,amstext}
\usepackage{dcolumn}
\usepackage{color}
\usepackage{subfigure}
\usepackage{graphicx}
\usepackage{epsfig,epstopdf}
\usepackage[applemac]{inputenc}  
\usepackage{pdfsync}
\usepackage{fancyhdr}

\usepackage{color}
\usepackage{xcolor}
\usepackage[breaklinks,colorlinks,backref]{hyperref}

\hypersetup{
    colorlinks, %set true if you want colored links
    linktoc=all, %set to all if you want both sections and subsections linked
    linkcolor=red,  %choose some color if you want links to stand out
  citecolor=hyptxt,
  urlcolor=blue}
  \hypersetup{
  citebordercolor=,
  filebordercolor=green,
  linkbordercolor=blue
}
\definecolor{hyptxt}{rgb}{0.7, 0.4, 0.9}

\usepackage{bibtopic}

\usepackage[full]{textcomp} % to get the right copyright, etc.
     \usepackage{fbb} % osf for text, lining for math
     \usepackage[varqu,varl]{inconsolata}% typewriter
     \usepackage[libertine,bigdelims,vvarbb]{newtxmath}
     \usepackage[cal=boondoxo]{mathalfa}% less slanted than STIX
    \usepackage[T1]{fontenc}
\useosf
\newtheorem{prop}{Proposition}[section]

\newcommand{\beprop}{\begin{prop}}
\newcommand{\enprop}{\end{prop}}
\newcommand{\bprf}{\begin{proof}}
\newcommand{\eprf}{\end{proof}}

\newcommand{\ds}{\displaystyle}

\newcommand{\ket}[1]{|\kern.3ex#1\kern.3ex\rangle}
\newcommand{\bra}[1]{\langle\kern.3ex #1 \kern.3ex|}
\newcommand{\scalar}[2]{\langle\kern.3ex #1 \kern.3ex|\kern.3ex#2\kern.3ex\rangle}

\def\N{\mathbb{N}}
\def\C{\mathbb{C}}

\def\lg{\langle }
\def\rg{\rangle }

\def\ii{\mathrm{i}}
\def\ud{\mathrm{d}}

\definecolor{hervecolor}{rgb}{0.8,0,0.7}

\begin{document}

\title{Holomorphic Discrete Series of SU$(1,1)$:\\Orthogonality Relations, Character Formulas, and Multiplicities in Tensor Product Decompositions}

\author{Jean-Pierre Gazeau}
\email{gazeau@apc.in2p3.fr}
\affiliation{APC, UMR 7164,\\Universit\'e Paris Cit\'e, 75013 Paris, France}

\author{Mariano A. del Olmo}
\email{marianoantonio.olmo@uva.es}
\affiliation{Departamento de F\'{\i}sica Te\'orica, At\'omica y \'Optica and IMUVA -- Mathematical Research Institute,\\Universidad de Valladolid, Paseo Bel\'en 7, 47011 Valladolid, Spain}

\author{Hamed Pejhan}
\email{pejhan@math.bas.bg}
\affiliation{Institute of Mathematics and Informatics, Bulgarian Academy of Sciences,\\Acad. G. Bonchev Str. Bl. 8, 1113, Sofia, Bulgaria}

\date{\today}

%%%%%%%%%%%%%%%%%%%%%%%%%%%%%%%%%%%%%%%%%%%%%%
\begin{abstract}
The SU$(1,1)$ group plays a fundamental role in various areas of physics, including quantum mechanics, quantum optics, and representation theory. In this work we revisit the holomorphic discrete series representations of SU$(1,1)$, with a focus on orthogonality relations for matrix elements, character formulas of unitary irreducible representations (UIRs), and the decomposition of tensor products of these UIRs. Special attention is given to the structure of these decompositions and the associated multiplicities, which are essential for understanding composite systems and interactions within SU$(1,1)$ symmetry frameworks. These findings offer deeper insights into the mathematical foundations of SU$(1,1)$ representations and their significance in theoretical physics.  
%\begin{description}
%\item[PACS numbers] 02.40.Tt, 02.30.Tb, 02.40.Tt, 03.65.Ca,03.65.Fd,
%\end{description}
\end{abstract}

\maketitle

%%%%%%%%%%%%%%%%%%%%%%%%%%%%%%%%%%%%%%%%%%%%%%
\section{Introduction}\label{intro}
The SU$(1,1)$ symmetry provides a fundamental mathematical framework for various applications in modern physics, particularly in quantum mechanics \cite{MORSE}, as well as in the study of coherent \cite{4, QO, 5,6} and squeezed states \cite{7} in quantum optics, quantum computation \cite{8}, and quantum cosmology \cite{9}. Moreover, SU$(1,1)$ symmetries naturally emerge in the description of symmetry structures in curved spacetimes, especially in $(1+1)$-dimensional Anti-de Sitter spacetime \cite{bacry, Barut1986, GazOlmo20, dSBook}, where they play a key role in conformal field theories \cite{Maldacena1999, CFT1} and quantum gravity scenarios \cite{Witten1998, QG}.  

As the simplest non-abelian, noncompact Lie group with a simple Lie algebra, SU$(1,1)$ admits three fundamental classes of UIRs: the principal, complementary, and discrete series \cite{Bargmann1947, vilenkin68}. Among these, the discrete series representations are particularly significant, as they play a crucial role in modeling quantum states \cite{perel86} and in characterizing symmetries in curved spacetimes \cite{bacry, Barut1986, GazOlmo20, dSBook}.  

In this work, we focus on the (holomorphic) discrete series representations of SU$(1,1)$, delving into the orthogonality properties of their matrix elements, character formulas, and--crucially--the decomposition of tensor products. A thorough understanding of the tensor product structure is essential for unraveling the interplay between SU$(1,1)$ symmetry and multipartite quantum systems, as well as for constructing physically relevant states in quantum optics and representation theory. In particular, by analyzing the multiplicities that emerge in these decompositions, we shed light on the internal structure of composite systems governed by SU$(1,1)$ symmetry, laying a solid foundation for future research in areas where this symmetry plays a pivotal role--from quantum information science to black hole physics and quantum cosmology.

%%%%%%%%%%%%%%%%%%%%%%%%%%%%%%%%%%%%%%%%%%%%%%
\section{The SU$(1,1)$ group and its representation in discrete series}
The SU$(1,1)$ group, understood as the double covering group of the kinematical group associated with the $(1+1)$-dimensional Anti-de Sitter spacetime \cite{bacry, Barut1986, dSBook}, is given by \cite{vilenkin68}:
\begin{equation} \label{groupSU11}
\mathrm{SU}(1,1) = \left\{ g = \begin{pmatrix}
\alpha & \beta \\
\bar\beta & \bar\alpha
\end{pmatrix}  \;;\; \alpha, \beta \in \C\, , \ \det g= \vert \alpha\vert^2 - \vert \beta\vert^2 = 1\right\}\,.
\end{equation}
In our notation, the ``bar'' symbol over an entity represents its complex conjugate. Any element \( g \in \mathrm{SU}(1,1) \) can be factorized as follows \cite{vilenkin68}:
\begin{align}
g &= \begin{pmatrix}
   \alpha & \beta \\
  \bar{\beta} & \bar{\alpha}
\end{pmatrix}= \begin{pmatrix}
e^{\ii \varphi/2} &  0  \\
0 &  e^{-\ii \varphi/2} 
\end{pmatrix}\begin{pmatrix}
  \cosh\tau/2   &  \sinh\tau/2 \\
  \sinh\tau/2   &  \cosh\tau/2
\end{pmatrix} \begin{pmatrix}
e^{\ii \psi/2}  &  0  \\
0  &  e^{-\ii \psi/2} 
\end{pmatrix}\,, 
\end{align}
where \( 0 \leq \varphi \leq 2\pi \), \(-2\pi \leq \psi < 2\pi\), and \(0 \leq \tau < \infty \). It provides the following parametrization for $\alpha$ and $\beta$:
\begin{equation} \label{coortau}
\alpha = \cosh\tau/2\; e^{\ii (\varphi+\psi)/2}\,, \quad \beta = \sinh\tau/2\; e^{\ii (\varphi-\psi)/2}\,. 
\end{equation}
In terms of the coordinates $\tau, \varphi, \psi$, the Haar measure on the unimodular group SU\((1,1)\) takes the form:
\begin{equation} \label{HaarSU11B}
\ud_{\texttt{Haar}} (g) = \frac{1}{8\pi^2}\,\sinh\tau\,\ud\tau\,\ud\varphi\,\ud\psi\,.
\end{equation}

On the other hand, the holomorphic discrete series UIRs of the SU$(1,1)$ group are realized, in a broad sense, on Fock-Bargmann Hilbert spaces of holomorphic functions defined on the unit disk:
\begin{equation}
\mathcal{D} \equiv \left\{z \in \C\,, \, \vert z \vert < 1 \right\} \sim \mathrm{SU}(1,1)/\mathrm{U}(1)\,,
\end{equation}
where the group U$(1)$ is the maximal compact subgroup of SU$(1,1)$, and it is realized here as matrices of the form:
\begin{align}\label{maxisub}
h(\theta) = \begin{pmatrix} e^{i\theta/2} & 0 \\ 0 & e^{-i\theta/2} \end{pmatrix}\,,\quad  0 \leq \theta < 4\pi\,.
\end{align}
It is important to note that \(\mathcal{D}\) is one of the four two-dimensional K{\"a}hler manifolds \cite{doubrovine82,perel86}.

To provide a more precise visualization of these representations, let \(\mathcal{FB}_{\eta}\)  (for a given \(\eta > 1/2\)) be the Fock-Bargmann Hilbert space of all analytic functions \(f(z)\) on \(\mathcal{D}\), that are square-integrable with respect to the following scalar product:
\begin{equation} \label{FBsu11inprod}
\lg f_1|f_2\rg = \frac{2 \eta -1}{2 \pi} \, \int_{\mathcal{D}} \overline{f_1(z)}\, f_2(z)\, (1-\vert z \vert^2)^{2 \eta -2}\, \ud^2z\,.
\end{equation}
An orthonormal basis in this context corresponds to the powers of \(z \), appropriately normalized:
\begin{equation} \label{orthbasdisk2}
e_{n} (z) \equiv \, \sqrt{\frac{(2 \eta)_n}{n!}}\, z^n\,, \quad \mathrm{with}\;  z\in\mathcal{D}\,,\;n \in \N\,,
\end{equation}
where $(2 \eta)_n\equiv \Gamma(2 \eta +n)/\Gamma(2 \eta)$ is the Pochhammer symbol. For a given $\eta = 1,\, 3/2,\, 2,\, 5/2,\, \dotsc$, the SU$(1,1)$ UIR, $U^{\eta}(g)$, acts on $\mathcal{FB}_{\eta}$ as:
\begin{equation} \label{UIRDSgsu11}
\mathcal{FB}_{\eta} \ni f(z) \;\mapsto\; \left(U^{\eta}(g)\, f\right)(z) = (-\bar\beta \, z + \alpha)^{-2\eta}\; f \left( \frac{\bar\alpha z - \beta}{-\bar\beta z + \alpha}\right)\,,
\end{equation}
with $g\in SU(1,1,)$ given by \eqref{groupSU11}. A countable collection of these UIRs forms what is known as the ``almost complete'' holomorphic discrete series for SU\((1,1)\) \cite{vilenkin68}. The term ``almost complete'' reflects the fact that the lowest representation, corresponding to \(\eta = 1/2\), requires special handling because the inner product defined in (\ref{FBsu11inprod}) does not exist in this case. Moreover, if we extended our consideration to the continuous range \(\eta \in [1/2, +\infty)\), we would be compelled to work with the universal covering group of SU\((1,1)\).

As detailed in Refs. \cite{miller68, GazOlmo20}, the matrix elements of the operator \( U^{\eta}(g) \) with respect to the orthonormal basis (\ref{orthbasdisk2}) are expressed in terms of Jacobi polynomials:
\begin{align} \label{mateldiscsu11Jac}
\ds U^{\eta}_{nn^{\prime}}(g) \equiv \lg e_n | U^{\eta}(g) | e_{n^{\prime}}\rg =&\, \left( \frac{n_<!\, \Gamma(2 \eta + n_>) }{n_>!\, \Gamma(2 \eta + n_<)} \right)^{1/2}\, \alpha^{-2 \eta - n_>}\, \bar{\alpha}^{n_<} \nonumber\\ \times&\, \left(\gamma(\beta,\bar \beta)\right)^{n_> -n_<}\; P^{(n_>-n_<\,,\, 2 \eta -1)}_{n_<}\left( 1-2\vert z\vert^2 \right)\,,
\end{align}
where:
\begin{equation}
z = \beta \bar{\alpha}^{-1} \in {\mathcal{D}}\,, \quad \gamma(\beta,\bar\beta) = \left\lbrace \begin{array}{cc}
- \beta  & \quad  \mbox{if}\;\; n_> =  n^{\prime} \\
\bar\beta & \quad  \mbox{if}\;\;  n_> = n
\end{array}\right.\,, \quad \mbox{and} \quad 
n_{\substack{
> \\
<}} = \left\lbrace \begin{array}{c}
\max \\
\min
\end{array}\right.\, \;\; (n,n^{\prime}) \geq 0\,.
\end{equation}
In terms of the parameters \(\tau, \varphi, \psi\) given in Eq. \eqref{coortau}, the above expression \eqref{mateldiscsu11Jac} takes the form:
\begin{align} \label{mateldiscsu11JacC}
\nonumber U^{\eta}_{nn^{\prime}}(g) &= 2^{\frac{n_<-n_>}{2}-\eta}\left( \frac{n_<!\, \Gamma(2 \eta + n_>) }{n_>!\, \Gamma(2 \eta + n_<)} \right)^{1/2}
(1-x)^{\frac{n_>-n_<}{2}}(1+x)^{\eta}\times\\
&\times P^{(n_>-n_<\, ,\, 2 \eta -1)}_{n_<}\left(x \right)\times\left\lbrace\begin{array}{cc}
(-1)^{n^{\prime}-n} e^{-\ii(\eta+n)\varphi}  e^{-\ii(\eta+n^{\prime})\psi} & \quad \mbox{if}\;\; n^{\prime}\geq n  \\
e^{-\ii(\eta+n^{\prime})\varphi}  e^{-\ii(\eta+n)\psi} & \quad \mbox{if}\;\; n\geq n^{\prime}    
\end{array}\right.\,,
\end{align}
where, for the sake of simplicity, we have defined $x = 1-2\tanh^2\tau/2\in [-1,1]$, and hence, $\tanh\tau/2= \vert z\vert= \sqrt{\frac{1-x}{2}}$.

The trace of the operator $U^{\eta}(g)$, for an arbitrary $g\in \mathrm{SU}(1,1)$ \eqref{groupSU11}, is given by (see Ref. \cite{GazOlmo20} and references therein):
\begin{equation}
\label{traceUng}
\chi^{\eta}(g)\equiv  \mathrm{tr}\,\left(U^{\eta}(g)\right)= \frac{1}{2}\,\left((\Re\alpha)^2-1\right)^{-1/2}\, \left(\Re\alpha + \left((\Re\alpha)^2-1\right)^{1/2}\right)^{1-2\eta}\, .  
\end{equation}
In terms of the coordinates $(x,\varphi,\psi)$, it reads as:
\begin{equation} \label{tracexpp}
\chi^{\eta}(g) = \frac{1}{2}(1+x)^{\eta} \left(\cos(\varphi+\psi)-x\right)^{-1/2}\,\left(\sqrt{2}\cos\left(\frac{\varphi+\psi}{2}\right)+ (\cos(\varphi+\psi)-x)^{1/2}\right)^{1-2\eta}\,. 
\end{equation}

%%%%%%%%%%%%%%%%%%%%%%%%%%%%%%%%%%%%%%%%%%%%%%
\section{Orthogonality relations} \label{orthogonalidad}
A fundamental property of the discrete series of UIRs of Lie groups, as expected from their nature, is their orthogonality relations. Specifically, for all pairs \((\eta_1, \eta_2)\), where \(2\eta_i \in \mathbb{N}\) and \(2\eta_i > 1\) for \(i = 1,2\), this orthogonality is a defining characteristic:
\begin{equation} \label{orthogrel2}
\int_{\mathrm{SU}(1,1)}\ud_{\mathrm{Haar}}(g)\, U^{\eta_1}_{mm^{\prime}}(g)\, \overline{U^{\eta_2}_{nn^{\prime}}(g)} = d_{\eta_1}\,\delta_{\eta_1 \eta_2}\, \delta_{mn}\, \delta_{m^{\prime}n^{\prime}} \,, 
\end{equation}
where $d_{\eta} = 2/(2\eta-1)$ is the (formal) dimension of the representation $U^{\eta}$. 

Let us verify this result by selecting, without loss of generality, \( m^\prime \geq m \) and \( n^\prime \geq n \). Utilizing the expression for the matrix elements given in Eq. \eqref{mateldiscsu11JacC}, the left-hand side of the above equality takes the form:
\begin{align}
&\int_{\mathrm{SU}(1,1)}\ud_{\mathrm{Haar}}(g)\, U^{\eta_1}_{mm^{\prime}}(g)\, \overline{U^{\eta_2}_{nn^{\prime}}(g)}= \frac{(-1)^{m^{\prime}-m +n^{\prime}-n}}{2\pi^2}\,2^{\frac{m-m^{\prime}+ n-n^{\prime}}{2}-\eta_1-\eta_2}\left( \frac{m!\,n!\, \Gamma(2 \eta_1 + m^{\prime})\Gamma(2 \eta_2 + n^{\prime}) }{m^{\prime}!\,n^{\prime}!\,\Gamma(2 \eta_1 + m)\, \Gamma(2 \eta_2 + n)} \right)^{1/2} \nonumber\\
&\qquad\qquad\qquad\times \int_{-1}^1\ud x \,(1-x)^{\frac{m^{\prime}-m+n^{\prime}-n}{2}}(1+x)^{\eta_1 +\eta_2-2}\,P^{(m^{\prime}-m\, ,\, 2 \eta_1 -1)}_{m}\left(x \right)\,P^{(n^{\prime}-n\, ,\, 2 \eta_2 -1)}_{n}\left(x \right) \nonumber\\
&\qquad\qquad\qquad\times \int_0^{2\pi}\,\ud\varphi\, e^{-\ii(\eta_1-\eta_2+m-n)\varphi} 
\int_{-2\pi}^{2\pi}\ud\psi \, e^{-\ii(\eta_1-\eta_2+m^{\prime}-n^{\prime})\psi}  \nonumber \\
&\qquad\qquad = 2^{2 + {m - m^{\prime}} -\eta_1-\eta_2}\left( \frac{m!\,n!\, \Gamma(2 \eta_1 + m^{\prime})\Gamma(2 \eta_2 + n^{\prime}) }{m^{\prime}!\,n^{\prime}!\,\Gamma(2 \eta_1 + m)\, \Gamma(2 \eta_2 + n)} \right)^{1/2}\,\delta_{\eta_1 -\eta_2,n-m}\,\delta_{\eta_1 -\eta_2,n^{\prime}-m^{\prime}} \nonumber\\
&\qquad\qquad\qquad\times \int_{-1}^1\ud x \,(1-x)^{m^{\prime}-m}(1+x)^{\eta_1 +\eta_2-2}\, P^{(m^{\prime}-m\, ,\, 2 \eta_1 -1)}_{m}\left(x \right)\, P^{(m^{\prime}-m\, ,\, 2 \eta_2 -1)}_{n}\left(x \right)\,. \label{111}
\end{align}
We proceed by considering the above relation in two distinct cases. First, when \(\eta_1 = \eta_2 = \eta\), meaning we examine the orthogonality relations between the matrix elements of a single UIR, \( U^\eta \). Second, when \(\eta_1 \neq \eta_2\), we investigate in this case the orthogonality relations between the matrix elements of two different UIRs, \( U^{\eta_1} \) and \( U^{\eta_2} \).

In the specific case where \(\eta_1 = \eta_2 = \eta\), the structure, as dictated by the delta functions in Eq. \eqref{111}, ensures that nonzero values arise exclusively when \( n = m \) and \( n^\prime = m^\prime \). Consequently, Eq. \eqref{111} reduces to:
\begin{align}
&2^{2 + m - m^{\prime} - 2\eta} \left( \frac{m!\, \Gamma(2 \eta + m^{\prime})}{m^{\prime}!\,\Gamma(2 \eta + m)} \right)\, \int_{-1}^1\ud x \,(1-x)^{m^{\prime}-m}(1+x)^{2\eta-2}\, \left( P^{(m^{\prime}-m\, ,\, 2 \eta -1)}_{m}\left(x \right) \right)^2 = \frac{2}{2\eta -1 }\,.
\end{align}
Note that the above result is derived using the Gradshteyn-Ryzhik identity (7.391) from Ref. \cite{gradryz07}:
\begin{equation*} \label{jac12}
\int_{-1}^{+1}\ud x \, (1-x)^{a} \, (1 +x)^{b-1} \, \left(P^{(a\, ,\, b)}_{m}(x)\right)^2\\
=\frac{2^{a + b}}{b}\,\frac{\Gamma(a + m +1)\,\Gamma(b + m +1)}{m!\,\Gamma(a +b + m +1)}\,,
\end{equation*}
for $a > -1$ and $b > 0$. 

Now, we focus on the second case, where \(\eta_1 \neq \eta_2\). Let us set \(\eta_1 = \eta_2 + s\) with \(s \in \mathbb{N}\) and \(s \geq 1\). Given the constraint \(\eta_1 - \eta_2 = n - m\), it follows that \(n = m + s\). Defining \(a = m^\prime - m \in \mathbb{N}\) and \(2\eta_2 - 1 = b \in \mathbb{N}_*\), our goal is to prove that:
\begin{equation} \label{nulljac2}
\int_{-1}^1\ud x \,(1-x)^{a}(1+x)^{b + s-1}\,P^{(a\, ,\, b+2s)}_{m}\left(x \right)\,P^{(a\, ,\, b)}_{m+s}\left(x \right)= 0\,. 
\end{equation}
Since any polynomial \( p_n(x) \) of degree \( n \) is orthogonal to any Jacobi polynomial of degree greater than \( n \), we deduce that: 
\begin{align}
\forall\, r=0,1,\cdots, m+s-1 \;\;:\;\; \int_{-1}^1\ud x \,(1-x)^{a}(1+x)^{b }\,x^r\,P^{(a\, ,\, b)}_{m+s}\left(x \right)= 0 \,.  
\end{align}
Furthermore, we note that \((1+x)^{s-1}P^{(a\,,\, b+2s)}_{m}(x)\) is a polynomial of degree \(m+s-1\). Consequently, Eq. \eqref{nulljac2} holds true.

%%%%%%%%%%%%%%%%%%%%%%%%%%%%%%%%%%%%%%%%%%%%%%
\section{Reduction of tensor products of two SU$(1,1)$ UIRs  of the discrete series}
Let us consider the tensor product $U^{\eta_1}\otimes U^{\eta_2}$ of two UIRs of the discrete series, which means that $\eta_1, \eta_2 \in \N^\ast/2$, where $ \N^\ast = \N- \{0\}$. The decomposition of this tensor product yields the following expansion:
\begin{equation} \label{redTP1}
U^{\eta_1}\otimes U^{\eta_2} = \bigoplus_{\eta_3\in  \N^\ast/2} \;\mathfrak{m}(1,2,3) \,U^{\eta_3}\,,
\end{equation}
where $\mathfrak{m}(1,2,3)$ is the multiplicity of the appearance of the UIR $U^{\eta_3}$ in the expansion. 
An immediate result is provided by the expression  \eqref{traceUng} for the character $\chi^{\eta}(g)$. Specifically, by restricting \( g \) to the maximal compact subgroup, i.e., setting \( g = h(\theta) \in U(1)\) as given in Eq. \eqref{maxisub}, we obtain the following expression for the character:
\begin{equation} \label{charht}
\chi^{\eta} (h(\theta)) = \frac{1}{2\ii \sin\theta/2} e^{\ii(1-2\eta)\theta/2}\,, 
\end{equation}
where we made the sign choice $\sqrt{\cos^2\theta/2 -1}= \ii \sin\theta/2$. Starting from the product:
\begin{equation} \label{expredchar1}
\chi^{\eta_1}(h(\theta))\, \chi^{\eta_2}(h(\theta))=- \frac{1}{4 \sin^2\theta/2} e^{\ii({1-\eta_1 - \eta_2})\theta}\,,
\end{equation}
and applying the expansion formula:
\begin{equation} \label{expisin}
\frac{1}{\sin\theta/2} = 2\ii e^{-\ii\theta/2}\sum_{n\geq 0} e^{-\ii n\theta}\,, 
\end{equation}
to Eq. \eqref{expredchar1} yield the following decomposition: 
\begin{equation} \label{expredchar2}
\chi^{\eta_1}(h(\theta))\, \chi^{\eta_2}(h(\theta)) =  \frac{1}{2\ii\sin\theta/2} \sum_{n\geq 0}  e^{\ii(1-2(\eta_1+\eta_2 + n))\theta/2} = \sum_{n\geq 0}\chi^{\eta_1 + \eta_2 + n}(h(\theta)) \,.
\end{equation}
It yields the following values for the multiplicity:
\begin{equation} \label{multired}
\mathfrak{m}(1,2,3) = \delta_{\eta_3, \eta_1 + \eta_2 + n}\, ,\quad n\in \N\,. 
\end{equation} 
Notably, this result both confirms and extends classical findings in the literature, such as those by Repka and others \cite{Repka, Tomasini}, and aligns with the established framework of Clebsch-Gordan decompositions for SU$(1,1)$.

%%%%%%%%%%%%%%%%%%%%%%%%%%%%%%%%%%%%%%%%%%%%%%
\section*{Acknowledgments} 
Mariano A. del Olmo is supported by  MCIN with funding from the European Union NextGenerationEU (PRTRC17.I1), and also by  PID2023-149560NB-C21 financed by MICIU/AEI/10.13039/501100011033 of Spain. Hamed Pejhan is supported by the Bulgarian Ministry of Education and Science, Scientific Programme ``Enhancing the Research Capacity in Mathematical Sciences (PIKOM)", No. DO1-67/05.05.2022. Jean-Pierre Gazeau would like to thank the University of Valladolid for its hospitality. This article/publication is based upon work from COST Action CaLISTA CA21109 supported by COST (European Cooperation in Science and Technology).

%%%%%%%%%%%%%%%%%%%%%%%%%%%%%%%%%%%%%%%%%%%%%%

\end{document}